\documentclass[10pt, leqno]{amsart}

\usepackage{amssymb,amsfonts,amscd,amsbsy}
\usepackage[mathscr]{eucal}

\theoremstyle{plain}
\newtheorem{thm}{Theorem}[section]

\newtheorem{pro}[thm]{Proposition}
\newtheorem{cor}[thm]{Corollary}
\theoremstyle{definition}

\newtheorem{rem}[thm]{Remark}

\newtheorem{conj}[thm]{Conjecture}

\begin{document}

\title[Binary Quadratic Forms]{Representations of integers by certain positive definite binary quadratic forms}
\author[Ram Murty and Robert Osburn]{Ram Murty and Robert Osburn}

\address{Department of Mathematics $\&$ Statistics, Queen's University, Kingston, Ontario, Canada K7L 3N6}

\email{murty@mast.queensu.ca}
\email{osburnr@mast.queensu.ca}

\subjclass[2000] {Primary 11E25, 11E45}

\begin{abstract}
We prove part of a conjecture of Borwein and Choi concerning an
estimate on the square of the number of solutions to $n=x^2+Ny^2$
for a squarefree integer $N$.
\end{abstract}

\maketitle
\section{Introduction}

We consider the positive definite quadratic form $Q(x,y)=x^2+Ny^2$ for a squarefree integer $N$. Let $r_{2,N}(n)$ denote the number of solutions to $n=Q(x,y)$ (counting signs and order). In this note, we estimate

\begin{center}
$\displaystyle\sum_{n\le x} {r_{2,N}(n)}^2$.
\end{center}

A positive squarefree integer $N$ is called solvable if $x^2+Ny^2$
has one form per genus. Note that this means the class number of
the form class group of discriminant $-4N$ equals the number of
genera, $2^t$, where $t$ is the number of distinct prime factors
of $N$. Concerning $r_{2,N}(n)$, Borwein and Choi \cite{BC} proved
the following:

\begin{thm}
Let $N$ be a solvable squarefree integer. Let $x >1$ and $\epsilon >0$. We have

\begin{center}
$\displaystyle\sum_{n \le x} {r_{2,N}(n)}^2 = \frac{3}{N} \Big( \displaystyle\prod_{p|2N} \frac{2p}{p+1} \Big ) (x\log x + \alpha(N)x) + O(N^{\frac{1}{4} + \epsilon} x^{\frac{3}{4} + \epsilon})$
\end{center}

\noindent where the product is over all primes dividing $2N$ and

\begin{center}
$\alpha(N) = -1 + 2\gamma + \displaystyle\sum_{p|2N} \frac{\log p}{p+1} + \frac{2L^{\prime}(1,\chi_{-4N})}{L(1, \chi_{-4N})} - \frac{12}{{\pi}^2} \zeta^{\prime}(2)$
\end{center}

\noindent where $\gamma$ is the Euler-Mascheroni constant and $L(1, \chi_{-4N})$ is the L-function corresponding to the quadratic character mod $-4N$.

\end{thm}

Based on this result, Borwein and Choi posed the following:

\begin{conj} For any squarefree $N$,

\begin{center}
$\displaystyle\sum_{n \le x} {r_{2,N}(n)}^2 \sim \frac{3}{N} \Big( \displaystyle\prod_{p|2N} \frac{2p}{p+1} \Big ) x\log x$
\end{center}

\end{conj}

Our main result is the following.

\begin{thm}
Let $Q(x,y)=x^2+Ny^2$ for a squarefree integer $N$ with $-N
\not\equiv 1 \bmod 4$. Let $r_{2,N}(n)$ denote the number of
solutions to $n=Q(x,y)$ (counting signs and order). Then
\begin{center}
$\displaystyle\sum_{n \le x} {r_{2,N}(n)}^2 \sim \frac{3}{N} \Big(\displaystyle \prod_{p|2N} \frac{2p}{p+1} \Big ) x\log x$.
\end{center}
\end{thm}

\section{Preliminaries}

We first discuss two key estimates and a result of Kronecker on
genus characters. Then using Kronecker's result, we prove a
proposition relating genus characters to poles of the
Rankin-Selberg convolution of L-functions. The first estimate is a
recent result of K$\ddot{u}$hleitner and Nowak \cite{Nowak},
namely

\begin{thm}

Let $a(n)$ be an arithmetic function satisfying $a(n) \ll n^{\epsilon}$ for every $\epsilon > 0$, with a Dirichlet series
\begin{center}
$F(s) = \displaystyle\sum_{n=1}^{\infty} \frac{a(n)}{n^s} = \frac{(\zeta_{K}(s))^2}{(\zeta(2s))^{m_1}(\zeta_{K}(2s))^{m_2}} G(s)$
\end{center}

\noindent where $\Re(s)>1$ and $\zeta_{K}(s)$ is the Dedekind zeta function of some quadratic number field $K$, $G(s)$ is holomorphic and bounded in some half plane $\Re(s) \geq \theta$, $\theta < \frac{1}{2}$, and $m_1$, $m_2$ are nonnegative integers. Then for $x$ large,
\begin{center}
$\displaystyle\sum_{n \leq x} a(n) = Res_{s=1} \Big (
F(s)\frac{x^s}{s} \Big ) + O(x^{\frac{1}{2}} (\log x)^{3} (\log
\log x)^{m_1 + m_2})$
\end{center}
\begin{center}
$= Ax\log x + Bx + O(x^{\frac{1}{2}} (\log x)^{3} (\log \log
x)^{m_1 + m_2})$
\end{center}

\noindent where $A$ and $B$ are computable constants.
\end{thm}

For an arbitrary quadratic number field $K$ with discriminant
$d_K$, let $\mathcal{O}_{K}$ denote the ring of integers in $K$,
and $r_{K}(n)$ the number of integral ideals $\mathcal{I}$ in
$\mathcal{O}_{K}$ of norm $N(\mathcal{I})=n$. From (4.1) in
\cite{Nowak}, we have
\begin{center}
$\displaystyle\sum_{n=1}^{\infty} \frac{(r_{K}(n))^2}{n^s} =
\frac{(\zeta_{K}(s))^2}{\zeta(2s)} \displaystyle\prod_{p|d_K}
(1+p^{-s})^{-1}$.
\end{center}

Applying Theorem 2.1 with $m_{1}=1$ and $m_{2}=0$, we obtain

\begin{cor} For any quadratic field $K$ of discriminant $d_K$ and $x$ large,

\begin{center}
$\displaystyle\sum_{n \leq x} (r_{K}(n))^2 = A_{1}x\log x + B_{1}x
+ O(x^{\frac{1}{2}} (\log x)^{3} \log \log x)$,
\end{center}

\noindent with $A_{1} = \frac{6}{{\pi}^{2}} {L(1, \chi_{d_K})}^2
\displaystyle\prod_{p|d_K} \frac{p}{p+1}$ and $B_1 = A_{1}
\alpha(N)$ with $\alpha(N)$ as in Theorem 1.1.
\end{cor}

The second estimate is a classical result of Rankin \cite{Rankin} and Selberg \cite{Sel1} which estimates the size of Fourier coefficients of
a modular form. Specifically, if $f(z)=\displaystyle\sum_{n=1}^{\infty} a(n) e^{2{\pi}inz}$ is a nonzero cusp form of weight $k$ on $\Gamma_{0}(N)$, then
\begin{center}
$\displaystyle\sum_{n \leq x} |a(n)|^{2} = \alpha \langle f, f \rangle x^{k} + O(x^{k-\frac{2}{5}})$
\end{center}

\noindent where $\alpha > 0$ is an absolute constant and $\langle f, f \rangle$ is the Petersson scalar product. In particular, if $f$ is a cusp form of weight 1, then $\displaystyle\sum_{n \leq x} |a(n)|^{2} = O(x)$.
One can adapt their result to say the following. Given two cusp forms of weight $k$ on a suitable congruence subgroup of $\Gamma=SL_{2}(\mathbb Z)$, say $f(z)=\displaystyle\sum_{n=1}^{\infty} a(n) e^{2{\pi}inz}$ and $g(z)=\displaystyle\sum_{n=1}^{\infty} b(n) e^{2{\pi}inz}$, then
\begin{center}
$\displaystyle\sum_{n \leq x} a(n)\overline{b(n)} n^{1-k} = Ax + O(x^{\frac{3}{5}})$
\end{center}

\noindent where $A$ is a constant. In particular, if $f$ and $g$ are cusp forms of weight 1, then $\displaystyle\sum_{n \leq x} a(n)\overline{b(n)} = O(x)$.

We will also use a result of Kronecker on genus characters. Let us
first explain some terminology. Let $K=\mathbb Q(\sqrt{d})$ be a
quadratic field of discriminant $d_K$. $d_K$ is said to be a prime
discriminant if it only has one prime factor. Thus it must be of
the form: $-4$, $\pm 8$, $\pm p \equiv 1 \bmod 4$ for an odd prime
$p$. Every discriminant can be written uniquely as a product of
prime discriminants, say $d_K=P_1 \dots P_k$. Here $k$ denotes the
number of distinct prime factors of $d_K$. Thus $d_K$ can be
written as a product of two discriminants, say $d_K=D_1D_2$ in
$2^{k-1}$ distinct ways (excluding order). Now, for any such
decomposition we define a character $\chi_{D_{1},D_{2}}$ on ideals
by

\begin{center}
$\chi_{D_{1},D_{2}}(\frak{p}) =\left \{ \begin{array}{l}
\chi_{D_1}(N\frak{p}) \quad \mbox{if $\frak{p} \nmid D_1$ } \\
\chi_{D_2}(N\frak{p}) \quad \mbox{if $\frak{p} \nmid D_2$ }
\end{array}
\right. $ \\
\end{center}

\noindent where $\chi_{d}(n)$ is the Kronecker symbol. This is
well defined on prime ideals because $\chi_{D}(N\frak{a})=1$ if
$(\frak{a}, D)=1$. $\chi_{D_{1},D_{2}}$ extends to all fractional
ideals by multiplicativity. Hence we have

\begin{center}
$\chi_{D_{1},D_{2}}: I \to \{\pm1\}$
\end{center}

\noindent where $I$ is the group of nonzero fractional ideals of
$\mathcal{O}_K$. Thus $\chi_{D_{1},D_{2}}$ has order two, except
for the trivial character corresponding to $d_K=d_K \cdot 1=1
\cdot d_K$. Every such character $\chi_{D_{1},D_{2}}$ is called
the genus character of discriminant $d_K$. As these are different
for distinct factorizations of $d_K$ (into a product of two
discriminants), we have $2^{k-1}$ genus characters. Kronecker's
theorem (see Theorem 12.7 in \cite{iwaniec}) is as follows.

\begin{thm} The L-function of $K$ associated with the genus character $\chi_{D_{1},D_{2}}$ factors into the Dirichlet L-functions,
\begin{center}
$L(s, \chi_{D_{1},D_{2}}) = L(s, \chi_{D_1}) L(s, \chi_{D_2})$.
\end{center}
\end{thm}

Let $K=\mathbb Q(\sqrt{-N})$, $N$ squarefree, $I$ as above, and
$P$ the subgroup of $I$ of principal ideals. For a non-zero
integral ideal $\frak{m}$ of $\mathcal{O}_K$, define

\begin{center}
$I(\frak{m}) = \{\frak{a} \in I : (\frak{a},  \frak{m})=1 \}$
\end{center}

\begin{center}
$P(\frak{m}) = \{\langle a \rangle \in P : a \equiv 1 \bmod
\frak{m} \}$.
\end{center}

A group homomorphism $\chi: I_{\frak{m}} \to S^{1}$ is an ideal
class character if it is trivial on $P(\frak{m})$, i.e.

\begin{center}
$\chi(\langle a \rangle )=1$
\end{center}

\noindent for $a \equiv 1 \bmod \frak{m}$. Thus an ideal class
character is a character on the ray class group $I(\frak{m})
\diagup P(\frak{m})$. Taking the trivial modulus $\frak{m}=1$, we
obtain a character on the ideal class group of $K$. Note that for
$K=\mathbb Q(\sqrt{-N})$ a genus character is an ideal class
character of order at most two.

Let us also recall the notion of the Rankin-Selberg convolution of
two L-functions. For squarefree $N$, consider two ideal class
characters $\chi_1$, $\chi_2$ for $\mathbb Q(\sqrt{-N})$ and their
associated Hecke L-series

\begin{center}
$L(s, \chi_1)= \displaystyle \sum_{n=1}^{\infty}
\frac{\chi_1(n)}{n^s}$
\end{center}

\begin{center}
$L(s, \chi_2)= \displaystyle \sum_{n=1}^{\infty}
\frac{\chi_2(n)}{n^s}$
\end{center}

\noindent which converge absolutely in some right half-plane. We
form the convolution L-series by multiplying the coefficients,

\begin{center}
$L(s, \chi_1 \otimes \chi_2) = \displaystyle \sum_{n=1}^{\infty}
\frac{\chi_1(n) \chi_2(n)}{n^s}$.
\end{center}

The following result describes a relationship between genus
characters $\chi$ and the orders of poles of $L(s, \chi \otimes
\chi)$. Precisely,

\begin{pro} Let $\chi$ be an ideal class character of $\mathbb
Q(\sqrt{-N})$, $-N \not\equiv 1 \bmod 4$, and $L(s,\chi)$ the
associated Hecke L-series. Then $\chi$ is a genus character if and
only if $L(s, \chi \otimes \chi)$ has a double pole at $s=1$.
\end{pro}

\begin{proof}
Suppose $\chi_{D_{1},D_{2}}$ is a genus character of discriminant
$-4N$, and $L(s, \chi_{D_{1},D_{2}}) = \displaystyle
\sum_{n=1}^{\infty} \frac{b_{i}(n)}{n^s}$. By Theorem 2.3 and
Exercise 1.2.8 in \cite{murty} (see the solution), we have

\begin{center}

$\displaystyle \sum_{n=1}^{\infty} \frac{{b_{i}(n)}^2}{n^s} =
\frac{L(s, \chi_{D_1}^2) L(s, \chi_{D_2}^2) L(s, \chi_{D_1}
\chi_{D_2})^2}{L(2s, \chi_{D_1}^2 \chi_{D_2}^2)}$.

\end{center}

Note that

\begin{center}
$L(s, \chi_{D_1}^2)= \zeta(s) \cdot \displaystyle \prod_{p|D_1}
(1-p^{-s})$,
\end{center}

\begin{center}
$L(s, \chi_{D_2}^2) = \zeta(s) \cdot \displaystyle \prod_{p|D_2}
(1-p^{-s})$,
\end{center}

\begin{center}
$L(s, \chi_{D_1} \chi_{D_2})^2=L(s, \chi_{-4N})^2$,
\end{center}

and

\begin{center}
$L(2s, \chi_{D_1}^2 \chi_{D_2}^2) = \zeta(2s) \cdot \displaystyle
\prod_{p|D_{1}D_{2}} (1-p^{-2s})$.
\end{center}

We have

\begin{center}
$\displaystyle \sum_{n=1}^{\infty} \frac{{b_{i}(n)}^2}{n^s} =
\frac{{\zeta(s)}^{2} L(s, \chi_{-4N})^2}{\zeta(2s)} \prod_{p|2N}
(1+p^{-s})^{-1}$
\end{center}

\noindent and thus a double pole at $s=1$.

Conversely, let $\chi$ be an ideal class character of $K=\mathbb
Q(\sqrt{-N})$ and suppose $L(s, \chi \otimes \chi)$ has a double
pole at $s=1$. Now $\chi$ is an automorphic form on
$GL_{1}(\mathbb{A}_{K})$. By automorphic induction (see
\cite{ac}), $\chi$ is mapped to $\pi$, a cuspidal automorphic
representation of $GL_{2}(\mathbb{A}_{\mathbb Q})$. Note that
$\pi$ is reducible as, otherwise, $L(s, \pi \otimes \pi)$ has a
simple pole at $s=1$ (\cite{ac}, page 200). As $K$ is a quadratic
extension of $\mathbb Q$, we must have $\pi = \chi_1 + \chi_2$
where $\chi_i$ are Dirichlet characters. As $L(s, \chi)=L(s, \pi)$
(see \cite{ac}) and thus $L(s, \chi \otimes \chi)= L(s, \pi
\otimes \pi)$,

\begin{center}
$L(s, \pi \otimes \pi) = L(s, \chi \otimes \chi) = \displaystyle
\frac{L(s, \chi_1^2) L(s, \chi_2^2) L(s, \chi_1 \chi_2)^2}{L(2s,
\chi_1^2 \chi_2^2)}$.
\end{center}

\noindent Now $L(s, \chi \otimes \chi)$ has a double pole at $s=1$
if and only if either $\chi_1=\overline{\chi_2}$, $\chi_2^2 \neq
1$, and $\chi_1^2 \neq 1$ or $\chi_1^2=1$, $\chi_2^2=1$, and
$\chi_1 \chi_2 \neq 1$. The latter implies $\chi$ is a genus
character. We now need to show that the former also implies that
$\chi$ is a genus character. Note that

\begin{center}
$L(s, \chi)= \displaystyle \prod_{\frak p} \Big(
1-\frac{\chi(\frak p)}{N(\frak p)^{s}}\Big)^{-1}$
\end{center}

\noindent and

\begin{center}
$L(s, \chi_1 + \chi_2) = \displaystyle \prod_{p} \Big(
1-\frac{\chi_1(p)}{p^{s}} \Big)^{-1}  \displaystyle \prod_{p}
\Big( 1-\frac{\chi_2(p)}{p^{s}} \Big)^{-1}$.
\end{center}

\noindent As $L(s, \chi)=L(s, \pi)$ and $L(s, \pi)=L(s, \chi_1 +
\chi_2)$, we compare Euler factors to get

\begin{center}
$\chi_1(p) + \chi_2(p) =\left \{ \begin{array}{l}
0 \quad \mbox{if $p$ is inert in $K$} \\
\chi(\frak{p}) + \overline{\chi(\frak{p})} \quad \mbox{if $p$
splits in $K$. }
\end{array}
\right. $ \\
\end{center}

For $p$ inert in $K$, this yields $\chi_1(p)=-\chi_2(p)$ and so
$\overline{\chi_2(p)} = \chi_1(p)=-\chi_2(p)$ which implies
$\chi_2^{2} (p)=-1$ and so $\chi_2(p)=\pm i$. Now consider the
following equation whose sum sieves the inert primes

\begin{center}

$\displaystyle \frac{1}{2} \sum_{\substack{p \le x \\  \\ p
\hspace{.05in} \text{prime}}} \Big (1 - \Big(\frac{-4N}{p} \Big)
\Big ) \chi_2^{2} (p) = -\pi(x)$.
\end{center}

\noindent Here $\pi(x)$ is the number of primes between $1$ and
$x$. Thus

\begin{center}

$\displaystyle \frac{1}{2} \sum_{\substack{p \le x \\  \\ p
\hspace{.05in} \text{prime}}} \chi_2^{2} (p)  - \frac{1}{2}
\sum_{\substack{p \le x \\  \\ p \hspace{.05in} \text{prime}}}
\Big(\frac{-4N}{p} \Big)  \chi_2^{2} (p) = -\pi(x)$.

\end{center}

As $\chi_2^2 \neq 1$, we have by the prime ideal theorem,
$\displaystyle \sum_{p \le x} \chi_2^{2}(p) = o(\pi(x))$ and so

\begin{center}

$\displaystyle \sum_{p \le x} \Big( \frac{-4N}{p} \Big
)\chi_{2}^{2}(p) \sim \pi(x)$.

\end{center}

This implies $\Big(\frac{-4N}{p} \Big)\chi_2^{2}(p) =1$. If $p$
splits in $K$, then $\chi_2^{2}(p)=1$ and so $\chi_2(p)=\pm 1$. A
similar argument works for $\chi_1$ and so we also have
$\chi_1(p)=\pm 1$ if $p$ splits in $K$.

Again comparing the Euler factors in $L(s, \chi)$ and $L(s, \pi)$,
the values of $\chi(\frak{p})$ must coincide with the values of
$\chi_1(p)$ and $\chi_2(p)$, that is, $\chi(\frak{p})=\pm 1$. Now
$\chi(\frak{p})=\chi([\frak{p}])$ where $[\frak{p}]$ is the class
of $\frak{p}$ in the ideal class group of $K$. By the analog of
Dirichlet's theorem for ideal class characters, we know that in
each ideal class $\frak{C}$ there are infinitely many prime ideals
which split. Thus $\chi(\frak{C})=\pm 1$ and hence is of order 2.
This implies $\chi$ is a genus character.

\end{proof}

\begin{rem}

By Proposition 2.4, if $\chi$ is a non-genus character, then $L(s,
\chi \otimes \chi)$ has at most a simple pole at $s=1$.

\end{rem}

\section{Proof of Theorem 1.3}

\begin{proof}
As $-N \not\equiv 1 \bmod 4$, the discriminant of $K=\mathbb
Q(\sqrt{-N})$ is $-4N$. We also assume that $t$ is the number of
distinct prime factors of $N$ and so the discriminant $-4N$ has
$t+1$ distinct prime factors.

Given the quadratic form $Q(x,y)=x^2+Ny^2$, we consider the
associated Epstein zeta function (see \cite{ep}, \cite{kani},
\cite{sc}, or \cite{siegel})

\begin{center}
$\zeta_{Q}(s) = \displaystyle\sum_{x,y \neq 0} \frac{1}{(x^2+Ny^2)^s} = \displaystyle\sum_{n=1}^{\infty} \frac{r_{2,N}(n)}{n^s}$.
\end{center}

\noindent for $\Re(s)>1$. Now for $K=\mathbb Q(\sqrt{-N})$, we have Dedekind's zeta function

\begin{center}
$\zeta_{K}(s) = \displaystyle\sum_{\frak{a}} \frac{1}{{N(\frak{a})}^s} = \displaystyle\sum_{n=1}^{\infty} \frac{a_n}{n^s}$
\end{center}

\noindent where the sum is over all nonzero ideals $\frak{a}$ of $\mathcal{O}_{K}$. We now split up $\zeta_{K}(s)$, according to the classes $c_{i}$ of the ideal class group $C(K)$, into the partial zeta functions (see page 458 of \cite{neu})
\begin{center}
$\zeta_{c_i}(s) = \displaystyle\sum_{\frak{a} \in c_{i}} \frac{1}{{N(\frak{a})}^s}$
\end{center}

\noindent so that $\zeta_{K}(s)= \displaystyle\sum_{i=0}^{h-1}
\zeta_{c_i}(s)$ where $h$ is the class number of $K$. In our case
$K=\mathbb Q(\sqrt{-N})$ is an imaginary quadratic field and so by
\cite{cox} (Theorem 7.7, page 137), we may write

\begin{center}
$\zeta_{K}(s) = \displaystyle\sum_{i=0}^{h-1} \zeta_{Q_i}(s)$
\end{center}

\noindent where $Q_i$ is a class in the form class group. Note
that in this context, $Q(x,y)$ corresponds to the trivial class
$c_{0}$ in $C(K)$ and so $\zeta_{c_{0}}(s) = \zeta_{Q(x,y)}(s)$.
Now let $\chi$ be an ideal class character and consider the Hecke
L-function for $\chi$, namely
\begin{center}
$L(s, \chi)= \displaystyle \sum_{\frak{a}}
\frac{\chi(\frak{a})}{{N(\frak{a})}^s}$
\end{center}

\noindent where $\frak{a}$ again runs over all nonzero ideals of
$\mathcal{O}_{K}$. We may now rewrite the Hecke L-function as
\begin{center}
$L(s,\chi) = \displaystyle\sum_{i=0}^{h-1} \chi(c_i)
\zeta_{c_{i}}(s)$.
\end{center}

And so summing over all ideal class characters of $C(K)$, we have
\begin{center}
$\displaystyle\sum_{\chi} \overline{\chi}(c_0) L(s, \chi) =
\displaystyle\sum_{i=0}^{h-1} \zeta_{c_i}(s) \Big (
\displaystyle\sum_{\chi} \overline{\chi}(c_0) \chi(c_{i}) \Big )$.
\end{center}

\noindent The inner sum is nonzero precisely when $i=0$. As
$\overline{\chi}(c_0)=1$ we have $\zeta_{c_0}(s) = \frac{1}{h}
\displaystyle\sum_{\chi} L(s, \chi)$. Thus
\begin{center}
$\zeta_{c_0}(s) = \frac{1}{h} (L(s, \chi_{0}) + L(s, \chi_{1}) +
\dots + L(s, \chi_{h-1}))$.
\end{center}

\noindent As $\chi_{0}$ is the trivial character, $L(s, \chi_{0})
= \zeta_{K}(s)$. We now compare $n^{th}$ coefficients, yielding

\begin{center}
$r_{2,N}(n) = \frac{1}{h} ( a_{n} + b_{1}(n) + \dots + b_{h-1}(n) )$
\end{center}

\noindent where $a_{n}$ is the number of integral ideals of
$\mathcal{O}_{K}$ of norm $n$ and the $b_{i}$'s are coefficients
of weight 1 cusp forms (see the classical work of Hecke
\cite{Hecke1}, \cite{Hecke2} or \cite{bump}). From the modern
perspective, this is straightforward. Each $L(s, \chi_{i})$, $1
\leq i \leq h-1$, can be viewed as an automorphic L-function of
$GL_{1}(\mathbb{A}_{K})$ and by automorphic induction (see
\cite{ac}) they are essentially Mellin transforms of (holomorphic)
cusp forms, in the classical sense. We now have

\begin{center}
$\displaystyle \sum_{n\le x} {r_{2,N}(n)}^2 = \frac{1}{h^2} \Big(
\sum_{n\le x} {a_{n}}^2 + \sum_{\substack{i \\ n\le x }}
{b_{i}(n)}^2 + 2 \sum_{\substack{i \\ n\le x }} a_{n}b_{i}(n) +
\sum_{\substack{i \neq j \\ n\le x}} b_{i}(n)b_{j}(n) \Big)$.
\end{center}

By the Rankin-Selberg estimate, $2 \displaystyle \sum_{\substack{i
\\ n\le x }} a_{n}b_{i}(n)$, $\displaystyle \sum_{\substack{i \neq j
\\ n\le x}} b_{i}(n)b_{j}(n)$ are equal to $O(x)$. By Corollary
2.2,

\begin{center}
$\displaystyle \frac{1}{h^2} \sum_{n\le x} {a_{n}}^2 =
\frac{1}{h^2} \Big( {A_1} x\log x + {B_1}x + O(x^{\frac{1}{2}}
(\log x)^{3} \log \log x) \Big )$.
\end{center}

We now must estimate $\displaystyle \sum_{\substack{i \\ n\le x }}
{b_{i}(n)}^2$. Let us now assume that the first $2^{t} - 1$ terms
arise from L-functions associated to genus characters. By
Proposition 2.4 and Nowak's proof of Theorem 2.1 (which uses
Perron's formula and the residue theorem), we obtain

\begin{center}
$\displaystyle \sum_{n\le x} {b_{i}(n)}^2 = A_{1}x \log x + B_{1}x + O(x)$
\end{center}

\noindent with $A_{1}$ and $B_{1}$ as in Corollary 2.2. As this
estimate holds for each $i$ such that $1 \le i \le 2^{t}-1$, the
term $A_1 x \log x$ appears $2^{t}$ times in the estimate of
$\displaystyle\sum_{n \le x} {r_{2,N}(n)}^2$. By Remark 2.5, the
remaining terms $\displaystyle \sum_{n\le x} {b_{i}(n)}^2$ for
$2^{t}-1 < i \le h-1$ are all $O(x)$. Thus

\begin{center}
$\displaystyle\sum_{n \le x} {r_{2,N}(n)}^2 = \frac{1}{h^2} \Big [ \Big ( 2^{t} \frac{6}{{\pi}^{2}} {L(1, \chi_{-4N})}^2 \displaystyle\prod_{p|2N} \frac{p}{p+1} \Big ) x \log x + O(x) \Big ] + O(x)$.
\end{center}

By (4.11) in \cite{gross} (or equation (8), page 171 in
\cite{cohn}), we have $L(1, \chi_{-4N})=\frac{h\pi}{\sqrt{N}}$ and
so

\begin{center}
$\displaystyle\sum_{n \le x} {r_{2,N}(n)}^2 = \frac{3}{N} \Big(\displaystyle \prod_{p|2N} \frac{2p}{p+1} \Big ) x\log x + O(x)$.
\end{center}

\noindent The result then follows.

\end{proof}

\begin{rem} It should be possible to obtain the second term
in the asymptotic formula. By a careful application of the
Rankin-Selberg method, one should obtain an error term of the form
$O(x^{\theta})$ with $\theta <1$. The remaining case $-N \equiv 1
\bmod 4$ requires more subtle analysis due to the fact that for
$K=\mathbb Q(\sqrt{-N})$, $\mathbb Z[\sqrt{-N}]$ is not the
maximal order of $K$. It involves the study of L-series attached
to orders. Using the techniques in \cite{cp} and \cite{kani}, we
will take this and sharper error terms up in some detail in a
forthcoming paper.
\end{rem}

\end{document}